\documentclass{article}
  \usepackage[reqno, namelimits, sumlimits]{amsmath}
        \usepackage{amssymb, amsfonts}
        \usepackage{mathrsfs}
        \usepackage[margin=1.4in, paperwidth=8.5in, paperheight=11.0in]{geometry}

\newcommand{\be}{\begin{equation}}
\newcommand{\ee}{\end{equation}}
\newcommand{\bb}{\bigskip}
\newcommand{\la}{\label}

\renewcommand{\div}{\,{\rm div}\,}

\newcommand{\Ad}{{\rm Ad}\,}

\newcommand{\rf}[1]{(\ref{#1})}

\newcommand{\hl}{\mathfrak h}

\newcommand{\ve}{\varepsilon}
\newcommand{\Om}{\Omega}

\newcommand{\vf}{\ensuremath{\varphi}}

\newcommand{\vek}[2]{\left(\begin{array}{c} #1 \\ #2 \end{array}\right)}

\newcommand{\R}{{\mathbf R}}

\newcommand{\AAA}{\mathcal A}

\newcommand{\OO}{\mathcal O}

\newcommand{\HH}{\mathcal H}

\newcommand{\MM}{\mathcal M}

\newcommand{\XX}{\mathcal X}
\newcommand{\YY}{\mathcal Y}

\newcommand{\om}{\omega}

\newcommand{\sm}{\smallskip}

\newcommand{\homog}[1]{$m-$homogeneous}

\newcommand{\pul}{\frac12}

\def\XXint#1#2#3{{\setbox0=\hbox{$#1{#2#3}{\int}$ }
\vcenter{\hbox{$#2#3$ }}\kern-.6\wd0}}

\newcommand{\ed}{\end{document}}

\newcommand{\edoc}{\end{document}}

\newcommand{\gl}{\mathfrak g}
\newcommand{\gls}{\mathfrak g^*}
\usepackage{amsmath}
\numberwithin{equation}{section}
\newtheorem{theorem}{Theorem}[section]

\newtheorem{corollary}{Corollary}[section]
\newtheorem{proposition}{Proposition}[section]
\newtheorem{definition}{Definition}[section]

\begin{document}
\title{Dynamics of geodesic flows with random forcing on Lie groups with left-invariant metrics }
\author{W.~Hu
 and V.~\v Sver\'ak\\ \\
 University of Minnesota
}

\date{}
\maketitle
\begin{abstract}We consider stochastic perturbations of geodesic flow for left-invariant metrics on \hbox{fi\-ni\-te-di\-men\-sional} Lie groups and study the H\"ormander condition  and some properties of the solutions of the corresponding Fokker-Planck equations.
\end{abstract}

\section{Introduction}
Our motivation for this paper comes from the problem of turbulent mixing. However, instead of studying the motion of fluids, which can be mathematically described by trajectories in the group of diffeomorphisms of the domain containing the fluid (as pointed out by V.~I.~Arnold~\cite{Arnold}), we will study its finite-dimensional version when the diffeomorphism group is replaced by a finite-dimensional  Lie group $G$. We equip $G$ with a left-invariant metric, and consider stochastically perturbed geodesic flows. This situation may admittedly be somewhat removed from important phenomena in the real flows related to very high (or even infinite) dimensionality of the phase-spaces relevant there (at least if do not take the dimension of $G$ as large parameter), but it does retain an important feature: the amplification of the stochastic effects by the non-linearity. A well-known example of this effect in the context of 2d incompressible fluids has been established in a seminal paper by Hairer-Mattingly~\cite{HM}, where ergodicity for stochastically forced 2d Navier-Stokes equation was proved for degenerate forcing, under optimal assumptions. There are at least two important themes involved in this result. One  might perhaps be called algebraic, and involves calculations of Lie algebra hulls related to the H\"ormander hypoellipticity condition~\cite{Hormander, Hairer-Hormander}. The other one belongs to Analysis/Probability and deals with consequences of the H\"ormader condition (which are of course already of great interest in finite dimension) in the infinite-dimensional setting (under suitable assumptions). In the finite-dimensional models we consider in this paper the analysis component is much simpler, although there still are many non-trivial and interesting issues  related to various aspects of hypoelliptic operators, such as the domains of positivity of the of fundamental solutions and convergence to equilibria.

Our focus here will be on the algebraic part. Roughly speaking, we will be interested in algebraic conditions which imply the H\"ormander condition, ergodicity and convergence to equilibria. The stochastic forces will be essential for this, but it is interesting to try in some sense to minimize the ``amount" of stochasticity which is needed.

One can of course also study the ergodicity of the non-stochastic dynamics, but we have nothing new to say about this notoriously difficult problem.

We consider two different types of models. The first one might be called the Langevin-type perturbation of the geodesic flow. It is related to the stochastic equation
\be\la{i1}
\ddot x +\nu \dot x = \xi\,,
\ee
where $\xi$ is a random force, which is ``degenerate", in the sense that it acts only in a few directions. On a group $G$ with a left-invariant metric (and under suitable assumptions on~$\xi$) one can employ symplectic reduction and obtain an equation
\be\la{i2}
\dot z^k = q^k(z,z)-\nu z^k + \sigma^k_l \dot w^l\,,
\ee
  in the Lie algebra $\gl$ of the group, where we sum over repeated indices, $k$ runs from $1$ to the dimension of the group, $l$ runs from $1$ to the dimension of the noise (which can be $1$), $w^l$ are independent standard Wiener processes, and the equation
 \be\la{i3}
 \dot z = q(z,z)
 \ee
 is the Euler-Arnold equation in $\gl $ as established in~\cite{Arnold}.
 For this model we determine an algebraic condition on $q$ which is necessary and sufficient for the H\"ormander condition for the corresponding Fokker-Planck equation to be satisfied in the cotangent space $T^*G$, see Theorem~\ref{thm1}. For a compact group $G$ this condition  implies ergodicity, and the projection of the ergodic measure to $G$ is the Haar measure. This means that the (stochastically perturbed) geodesic flow will visit all points on the group with the same probability (with respect to the Haar measure). We note that in the setting of the left-invariant metrics on a group this can never be the case without forcing, due to known conserved quantities one gets from Noether's theorem.

 For our next group of models we take a compact manifold $Z\subset \gl$ which is invariant under the flow of~\rf{i3} and consider
 \be\la{i4}
\dot z = q(z,z)+\xi\,,
\ee
where $\xi$ schematically stands for random forcing induced by the Brownian motion in $Z$ with respect to a natural Riemannian metric. One example we have in mind - in the co-tangent bundle $T^*G$ picture\footnote{We can of course go back and forth between $TG$ and $T^*G$ with the help of the metric.} -  is the intersection of a co-adjoint orbit and an energy level. The manifold $Z$ can have much lower dimension than $G$. This situation may in fact be a fairly realistic description  of a motion with random perturbations in which the quantities defining $Z$ are monitored and kept close to constant values by some control mechanism. When combined with  random perturbations, such control might easily induce random drift along the surface defined by specified values of the controlled quantities. (A more concrete mathematical process is described in subsection~\ref{constr33}.)
 Together with the equation
\be\la{kinematic}
a^{-1}\dot a =z,
 \ee
 the stochastic equation~\rf{i4} gives a stochastic equation in $G\times Z$.  In this situation we again determine an algebraic condition on $Z$ which is necessary and sufficient for the Fokker-Planck equation in $G\times Z$ associated to~\rf{i4} to satisfy the H\"ormander condition, see Theorem~\ref{thm2}, and, when the condition is satisfied, establish ergodicity and convergence to equilibrium. For compact $G$ the ergodic measure is given by a product of the Haar measure on $G$ and an invariant measure on $Z$.

  In the case of a non-compact $G$ and a compactly supported initial condition for the Fokker-Planck equation the behavior will of course be different,  and we illustrate what one might expect by an explicit calculation for $G=\R^n$ and a one-dimensional manifold $Z$, see Proposition~\ref{prop1}.

  The themes above have strong connections to control theory. In addition to the remark above about intrepreting $Z$ as a ``control surface", there is another connection via the Stroock-Varadhan Theorem~\cite{StroockVaradhan}. Roughly speaking, instead of random forcing $\xi$ one can consider forcing by control and ask which states can be reached (and how efficiently). For an introduction to control theory see for example~\cite{Jurdjevic}.

\section{Preliminaries}
\subsection{Basic notation and setup}
Let $G$ be a Lie group. Its elements will be denoted by $a, b, \dots$ We will denote by $\gl$ and $\gls$  respectively its Lie algebra and its dual.  Let $e_1,\dots, e_n$ be a basis of $\gl$ and let $e^1,\dots,e^n$ be its dual basis in $\gls$, determined by $\langle e^i,e_j\rangle=\delta^i_j$. We assume that a metric tensor with coordinates $g_{ij}$ in our basis is given on $\gl$. In what follows we will mostly use the standard formalism of othonormal frames and assume that $g_{ij}=\delta_{ij}$, which can of course always be achieved by a suitable choice of the original basis, although sometimes it may be  useful not to normalize $g_{ij}$ this way, so that other objects could be normalized instead. When $g_{ij}=\delta_{ij}$, we can then identify vectors with co-vectors without too much notation and write $|x|^2$ for the square of the norm  of an $x\in\gl$ or $x\in \gls$ given by the metric tensor.
However, we will try to avoid relying on this normalization too much and many of our formulae will be independent of it. In such situations we will use the classical convention of using upper indices for vectors and lower indices for co-vectors, with the usual convetions $y_k=g_{kl}y^l$ and $y^k=g^{kl}y_l$, where $g^{kl}$ is the inverse matrix of $g_{kl}$. In this notation we can for example write $|y|^2=y_ky^k$.

The various objects on $\gl$ and $\gls$ can be transported to $T_aG$ and $T_a^*G$ for any $a\in G$ in the standard way, by using the left translation $b\to ab$. The resulting frame of vectors fields on $G$ (or 1-forms), will still be denoted by $e_1,\dots, e_n$.

We can then consider $G$ as a Riemannian manifold. The left translations $b\to ab$ are more or less by definition isometries of the manifold. They obviously act transitively on $G$, and hence $G$ is a homogeneous Riemannian manifold.

The relevance of this construction for the mechanics of fluids and rigid bodies was pointed out in Arnold's paper~\cite{Arnold} already mentioned above.  The main point is that for fluids and rigid bodies the configuration space of the corresponding physical system is naturally given by a group (which however is infinite-dimensional for fluids), and the kinetic energy given a natural metric tensor on it. We  refer the reader to the book by Arnold and Khesin~\cite{ArnoldKhesin} for a deeper exposition of these topics and additional developments.

\subsection{The symplectic structure in $T^*G$ in left-invariant frames}
The cotangent space $T^*G$ is the canonical phase space for describing the geodesic flow in $G$ via the Hamiltonian formalism. For a group $G$ with a left-invariant metric the space $T^*G$ can be identified with $G\times \gls$ by using the frame $e_1,\dots,e_n$ on $G$:
\be\la{1}
(a,y)\in G\times \gls \qquad\rightarrow\qquad y_k e^k(a)\in T^*_aG\,,
\ee
where $e^1,\dots, e^n$ is the frame in $T^*G$ which is dual to $e_1,\dots,e_n$.
Here and in what follows we use the standard convention of summing over repeated indices. The ``coordinates" in $T^*G$ given by $(a,y)$ are convenient  for calculations, and will be freely used in what follows. Note that the prolongation of the action $a\to ba$ of $G$ on itself to $T^*G$ has a very simple form in the $(a,y)$ coordinates:
\be\la{1b}
(a,y)\to (ba,y)\,,
\ee
i.\ e.\ the $y$ coordinate stays unchanged. This is exactly because the frame $e^k$ is left-invariant.

As any cotangent space of a smooth manifold, the space $T^*M$ carries a natural symplectic structure. We start with the canonical 1-form on $T^*G$, which is given by
\be\la{2}
\alpha=y_k e^k(a)\,.
\ee
The symplectic form $\om$ is then given by
\be\la{3}
\om=d\alpha\,.
\ee
We have
\be\la{4}
d\alpha=dy_k\wedge e^k+y_kde^k\,.
\ee
The calculation of $de^k$ is standard. First, we introduce the structure constants of $\gl$ (with respect to the basis $e_k$ by
\be\la{5}
[e_i,e_j]=c^k_{ij}e_k\,.
\ee
Next, we apply Cartan's formula for the exterior differentiation:
\be\la{cartan}
de^k(e_i,e_j)=e_i\cdot e^k(e_j)-e_j\cdot e^k(e_i)-e^k([e_i,e_j])\,.
\ee
Combining~\rf{5} and~\rf{cartan}, together with the fact that the first two terms on the right-hand side of~\rf{cartan} vanish due to left-invariance of the objects involved, we obtain
\be\la{8}
\om\,\,=\,\,d\alpha\,\,=\,\,dy_k\wedge e^k-\pul \,y_k\,c^k_{ij}\,e^i\wedge e^j\,.
\ee
In other words, in the local frame on $T^*G$ given by $e_1,\dots, e_n,e^1\sim \frac{\partial}{\partial y_1},\dots, e^n\sim\frac{\partial}{\partial y_n}$, the form $\omega$ is given by the block matrix
\footnote{We use the usual identifications: if $f, g$ are two co-vectors with coordinates $f_i, g_j$ respectively, then the two-form $f\wedge g$ is identified with
the antisymmetric matrix $\om_{ij}=f_ig_j-f_jg_i$ and $(f\wedge g)(\xi,\eta)=
\om_{ij}\xi^i\eta^j$ for any two vectors $\xi,\eta$.}
\be\la{9}
\left(\begin{array}{cc}
 -C(y) & -I \\ I & 0
\end{array}\right)\,,
\ee
where $C(y)$ denotes the matrix $y_kc^k_{ij}$. The inverse of the matrix~\rf{9} is
\be\la{10}
\left(\begin{array}{cc}
 0 & I \\ -I & -C(y)
\end{array}\right)\,,
\ee
and for any function $H=H(a,y)$ on $T^*G$ the corresponding Hamiltonian equations are
\be\la{H}
\begin{array}{rcl}
(a^{-1}\dot a)^k & = & \frac {\partial H}{\partial y_k}\,,\\
\dot y_k & = & -e_k H +y_lc^l_{jk}\frac{\partial H}{\partial y_j}\,,
\end{array}
\ee
where $(a^{-1}\dot a)^k$ denotes the $k-$th coordinate of the vector $a^{-1}\dot a\in \gl$, the expression $e_k H$ denotes the derivative along the $e_k$ direction in the variable $a$. The last term on the right-hand side of the second equations represents the  Poisson bracket $\{H,y_k\}$ with $H$ considered as a functions of $y$  (and $a$  considered as fixed when calculating the bracket). The bracket is uniquely given by its usual properties and the relations
\be\la{PB}
\{y_i,y_j\}=y_kc^k_{ij}\,.
\ee
 It can be obtained by applying the standard Poisson bracket on the symplectic manifold $T^*G$ to functions independent of $a$ in the above coordinates $(a,y)$.

 Note that the equations~\rf{H} do not depend on the metric, they depend only on the structure of the Lie algebra.

\subsection{The symplectic reduction to $\gls$ and the Euler-Arnold equation}

When $H$ is invariant under the prolongation of the action by left multiplication of $G$ on itself to $T^*G$, which is equivalent to $H$ not depending on $a$ in the above coordinates $(a,y)$, i.\ e.\ $H=H(y)$, then the second equation of~\rf{H} does not contain $a$ and is simply
\be\la{EA}
\dot y_k=\{H,y_k\}\,.
\ee
This is a form of the Euler-Arnold equation, originally formulated in $\gl$ in~\cite{Arnold}. This equation represents one form of the reduction of the equations on $T^*G$ to $\gls$ by the symmetries of the left action of $G$ on itself, see for example~\cite{MarsdenWeinstein}. The space $\gls$ has a natural structure of a Poisson manifold (with the Poisson bracket given by~\rf{PB}) and is foliated into ``symplectic leaves", which are given by the orbits of the co-adjoint representation, see for example~\cite{ArnoldKhesin, MarsdenWeinstein}. The orbits are given by
\be\la{orb}
\OO _{\bar y}=\{(\Ad a)^* \bar y, \,\,a\in G\}
\ee
where $\bar y$ is a fixed vector in $\gls$ and $\Ad a$ is defined below, and they have a natural structure of a symplectic manifold (with the maps $(\Ad a)^*$ acting by as symplectic diffeomorphism).

\subsection{Conserved quantities, the moment map, and Noether's theorem}
The Killing fields associated with the  symmetries of the Riemannian structure on $G$ with the left-invariant metric given by left multiplications $b\to ab$ are easily seen to be given by {\it right-invariant} vector fields $e(a)=\xi a$ (where $\xi\in \gl$) on $G$. By Noether's theorem there should be a conserved quantity associated to any such field. It is easy to see that the quantity is given by
\be\la{moment}
(a,y)\to (({\rm Ad} \,a^{-1})\,\xi\, ,\, y)=(\xi\,,\,({\rm Ad}\, a^{-1})^* y)\,,
\ee
where the operator $\Ad a$ is defined as usual by
\be\la{ad}
\Ad a \,\xi = a\xi a^{-1}\,.
\ee
The map $M\colon T^*G\to \gls$ given in the $(a,y)$ coordinates by
\be\la{M}
M(a,y)= (\Ad a^{-1})^*\, y
\ee
is the usual moment map associated with the (symplectic) action of $G$ on $T^*G$ (given by the prolongation of the left multiplication). The vector $M(a,y)$ is conserved under the Hamiltonian evolution, and the quantities $(\xi, M)$ are the conserved quantities from Noether's theorem applied to our situation. In particular, the Hamiltonian equations~\rf{H} obtained from taking the Hamiltonian as
\be\la{NH}
H(a,y)=(\xi, M(a,y))
\ee
are
\be\la{NHE}
\begin{array}{rcl}
(\dot a)a^{-1} & = & \xi \\
\dot y & = & 0\,.
\end{array}
\ee
The conservation of $M$ also has a geometric interpretation: if $x(s)$ is a geodecics (pa\-ra\-metrized by length) on a Riemannian manifold and $X$ is a Killing field (infinitesimal symmetry), then the scalar product $(\dot x, X)$ is constant. This is of course just another way to state the Noether's theorem in this particular case, but it can also be interpreted in terms of properties of Jacobi fields along our geodesics.

In the context of rotating rigid bodies, the quantity $M$ corresponds to the conservation of angular momentum, see~\cite{Arnold}. In the context of ideal fluids, the conservation of $M$ corresponds to the Kelvin-Helmholtz laws for vorticity, as observed by many authors.

It is easy to check the following fact: when $H$ is independent of $a$, i.\ e.\ $H=H(y)$, then for a curve $(a(t), y(t))$ in $T^*G$ satisfying the ``kinematic" equation
\be\la{kinem}
(a^{-1}\dot a)^k=  \frac {\partial H}{\partial y_k}
\ee
the ``dynamical" equation
\be\la{dyn}
\dot y_k = \{H,y_k\}
\ee
is equivalent to the (generalized) momentum conservation
\be\la{mc}
M(a,y)={\rm const.}
\ee

Also, if $(a(t),y(t))$ is a solution of the equations of motion and $H=H(y)$, then $y(t)$ is given by
\be\la{yt}
y(t)=(\Ad a(t))^*\,\, \bar y
\ee
for some fixed co-vector $\bar y\in\gls$.
\section{Perturbations by random forces}
\subsection{Langevin equation}
A very natural random perturbation of the geodesic flow is the Langevin equation, which can be symbolically written as
\be\la{L1}
\ddot a = - \nu \dot a + \ve \dot w\,,
\ee
for some parameters $\nu>0$ and $\ve >0$, which for a given $t>0$ and $a(t)$ is considered as an equation in $T_{a(t)} G$, with $\ddot a$ interpreted as the covariant derivative of $\dot a$ along the curve $a(t)$, and $w$ is a suitable form of Brownian motion in the Riemannian manifold $G$. Of course, the expression $\dot w$ is somewhat ambiguous and there are some subtle points in writing things in the correct way from the point of view of rigorous stochastic calculus. In particular, one has to distiguish carefully between the It\^o and Stratonovich integrals.
 Here we will mostly avoid the subtleties of the right interpretation of the stochastic equations such as~\rf{L1} by working instead with the Fokker-Planck equation, and we can define the transition probabilities for our processes via that equation.

A good starting point for writing the Fokker-Planck equation associated with~\rf{L1} is the Liouville equation in $T^*G$. This equation describes the evolution of a density $f(a,y)$ (with respect to the volume element given by the natural extension of the Riemannian metric from $G$ to $T^*G$, which is proportional to the volume element given by the $n-$th power $\om\wedge\om\wedge\dots\wedge \om$ ($n$ times) of the canonical symplectic form $\om$ above. The Liouville equation is
\be\la{Liouville}
f_t+ v^ke_kf+b_k\frac{\partial f}{\partial y_k}=0\,,
\ee
where
\be\la{vb}
v^k=\frac{\partial H}{\partial y_k}\,,\qquad b_k=\{H,y_k\}\,,
\ee
and $e_kf$ denotes the differentiation of $f(a,y)$ as a function of $a$ in the direction of the field $e_k$ defined earlier.
The vector field $X=v^ke_k+b_k\frac{\partial}{\partial y_k}$ is div-free (with respect to our volume form in $T^*G$), as follows from the Liouville theorem in Hamiltonian mechanics.
Hence equation~\rf{Liouville} is the same as
\be\la{L2}
f_t+\div (X f)=0\,,
\ee
where $\div$ is taken in our metric on $T^*G$. Our Fokker-Planck equation should then be
\be\la{FP}
f_t+v^ke_kf+b_k\frac{\partial f}{\partial y_k}+\frac{\partial}{\partial y_k}\left(-\nu v_k f -\frac{\ve^2}{2} \frac{\partial f}{\partial y^k}\right)=0 \,,
\ee
where $v_k$ is as above.
It can be considered as a combination of the Liouville transport with an Ornstein-Uhlenbeck process along the linear fibers of $T^*G$.  To have an exact correspondence to~\rf{L1}, we have to take $H(y)=\pul|y|^2$, the Hamiltonian of the geodesic flow.
 This equation can then be interpreted as describing a ``physical Brownian motion" of a particle in $G$. (We can for example think of $G$ being filled with an imcompressible fluid which is at rest and the Brownian particle suspended in the fluid and being subject to random ``kicks" from the fluid molecules and friction due to viscosity of the fluid, in the spirit of Einstein's paper~\cite{Einstein1905}.

 The symmetry reduction for~\rf{FP} corresponding the the symmetry reduction for~\rf{L1} is very simple: we consider it only for functions depending on $y$, which results in dropping the term $v^ke_kf$. The symplectic reduction of~\rf{H} to~\rf{EA} corresponds to the same procedure applied to the Liouville equation~\rf{Liouville}.

There is an explicit steady solution of~\rf{FP} given by
\be\la{sol1}
f(a,y)=Ce^{-\beta H }\,,\qquad \beta=\frac{2\nu}{\ve^2}\,,
\ee
where $C$ is any constant. The formula is the same as in the flat space. The approach to equilibrium will however  be influenced by the term $b_k\frac{\partial }{\partial y_k}$ which is absent in the flat case. Strictly speaking, the last statement applies unambiguously only to compact groups $G$, where the equilibrium~\rf{sol1} is easily seen to be unique among probability densities (for a suitable $C$, under some natural assumptions on $H$). We will discuss this point in some detail below in the more difficult case of degenerate forcing.

Given that the conservation of $M(a,y)$, from the point of view of Statistical Mechanics it is natural to consider (at least when $G$ is compact) distributions
in the phase space $T^*G$ given by
\be\la{equilibria}
f(a,y)=C e^{-\beta H(y) + (\xi\,,\, M(a,y))}=Ce^{-\beta H(y) + ((\Ad a^{-1})\xi\,,\, y)}
\ee
for $\beta>0$ and $\xi\in\gl$. In fact, if we replaced the Langevin equation by the Boltzmann equation
\be\la{Boltzmann}
f_t+v^ke_kf_k+b_k\frac{\partial f}{\partial y_k}=Q(f,f) \,,
\ee
for appropriate ``collision operator" $Q$ (defined on each fiber $T_a^*G$ in the same way as in the flat case), densities~\rf{equilibria} should be among the equilibria
(the set of which could possibly be larger due to symmetries other than those generated by the left shifts).
The large degeneracy of the set of equilibria is an important feature of the Boltzmann equation which is crucial for fluid mechanics. It is not shared by the Langevin equation, for which the equilibrium is unique (under reasonable assumptions). This is related to the hypoellipticity of the differential operator in~\rf{FP}, which we will discuss in some detail for more general operators in the next subsection.

We remark that one can modify the Langevin equation and get~\rf{equilibria} as equilibria for the modified equation. For this we simply change the Hamiltonian in~\rf{FP} to the expression
\be\la{nH}
\tilde H(a,y)= H(y)-((\Ad a^{-1})\xi\,, \,y)
\ee
This corresponds to watching a Brownian motion of a particle in incompressible fluid which moves in $G$ as a rigid body along the Killing field $\xi a$. (This is a steady solution of the equations of motion for an incompressible fluid.) The term
$((\Ad a^{-1})\xi\,, \,y)$ in the Hamiltonian is then produces the analogues of centrifugal and Coriolis forces which we encounter in rotating coordinate frames.

\subsection{Langevin equation with degenerate forcing}\label{sslang}
In PDEs of fluid mechanics one sometimes considers forcing through low spatial Fourier modes which is ``white noise" in time. See for example~\cite{HM, Kuksin}. In our context here this is akin to considering the system
\be\la{DL1}
\begin{array}{rcl}
(a^{-1}\dot a)^k & = & \frac{\partial H}{\partial y_k}\,\,,\\
\dot y_k & = & -e_k H + y_lc^l_{jk}\frac{\partial H}{\partial y_j}-\nu \frac{\partial H}{\partial y_k} + \sum_{i=1}^r \ve\dot w_i \tilde f^i_k\,\,,
\end{array}
\ee
where $\tilde f^1,\dots \tilde f^r$ are some fixed vectors in $\gls$ and $w_i$ are standard independent Wiener processes.
The term $-\nu \frac{\partial H}{\partial y_k}$ represents friction. One could consider more general forms of friction, but here we will be content with the above  special form. Many of the results below hold for more general friction forces.
The main complication in~\rf{DL1} as compared to the previous section is that $r$ can be less than the dimension $n$ of $\gls$.

In the remainder of this subsection we will assume that
\be\la{HH}
H=H(y)=\pul |y|^2=\pul y_ky^k\,,
\ee
which correspond to geodesic flow, or kinetic energy in classical mechanics.
Also, below we will need to do some Lie bracket calculations for which some formulae seem to be easier when we work in $\gl$ rather than $\gls$. This amounts to ``raising the indices" in the old-fashioned language, i.\ e.\ working in the coordinates $y^k$ rather then $y_k$. We note that with these assumptions we have
\be\la{vey}
y^k=v^k\,.
\ee
Equation~\rf{DL1} then becomes

\be\la{DL2}
\begin{array}{rcl}
(a^{-1}\dot a)^k & = & y^k\,\,,\\
\dot y^k & = & \tilde q^k_{ij}y^iy^j-\nu y^k + \sum_{i=1}^r\ve \dot w^i f_i^k \,\,,\,
\end{array}
\ee
where the notation is self-explanatory, perhaps with the exception of the term
$\tilde q^k_{ij}y^iy^j$, in which the coefficients are not uniquely determined by the function $y\to \tilde q(y,y)$.
A~straightforward ``raising of indices" gives the definition
\be\la{qdef}
([x,y],z)=(\tilde q(z,x),y)\,,\qquad x,y,z\in \gl\,\,,
\ee
which coincides with the Arnold form $B$ from~\cite{Arnold}. In what follows it will be advantageous to work with the symmetrization of $\tilde q$, which will be denoted by $q$:
\be\la{symq}
q(x,y)=\pul(\tilde q(x,y) + \tilde q(y,x))\,.
\ee
In  equation~\rf{DL2} it does not matter whether we use $\tilde q$ or $q$, of course.
Instead of~\rf{DL2} we can write
\be\la{dl3}
\begin{array}{rcl}
\dot a & = &  az\,\,,\\
\dot z & = &q(z,z)-\nu z + \ve \sigma \dot w\,\,,
\end{array}
\ee
where we use $z$ to emphasize that the equations are considered in $\gl$, as the variable $y$ was use to denote elements of $\gls$, $w$ is the vector of the standard Wiener process in $R^r$ and $\sigma$ is a suitable $n\times r$ matrix.
The corresponding Fokker-Planck equation for $f=f(a,z;t)$ then is
\be\la{fpd}
f_t+z^ke_k f + q^k(z,z)\frac{\partial f}{\partial z^k} + \frac{\partial}{\partial z^k} \left(-\nu z^k f- \frac{\ve^2}2 h^{kl} \frac{\partial f}{\partial z^l} \right)\,\,,
\ee
for a suitable symmetric positive semi-definite matrix $h$ (which is constant in $z$).

This is clearly a degenerate parabolic operator and we will study the classical parabolic H\"ormander condition for hypoellipticity for the Lie brackets generated by the vector fields relevant for the operator, see~\cite{Hairer-Hormander}.
In our context here the condition can be formulated in terms of the ``Lie algebra hull" of the vector fields
\be\la{add1}
\XX_k=\sigma_k^l\frac{\partial}{\partial z^l}
\ee
which satisfies the crucial additional condition that it is closed under the operation
\be\la{add2}
\Ad \XX_0\colon \XX\to [\XX_0,\XX]\,,
\ee
where the last bracket is the Lie bracket of vector fields and $\XX_0$ is the ``Euler-Arnold component" of our operator as specified below.

 The coordinates on $TG$ we will use are $(a,z)$, which correspond to $z^ke_k(a)$. The vector fields on $TG$ which will be relevant for our purposes will be of the form $A^k(z)e_k(a) + X^k(z)\frac{\partial}{\partial z^k}$\,. We will write
\be\la{fields}
A^k(z)e_k(a) + X^k(z)\frac{\partial}{\partial z^k} = \vek A X = \vek {A(z)}{X(z)}\,\,.
\ee
In these coordinates the Lie bracket is
\be\la{lbc}
\left[\,\,\vek A X \,\,,\,\,\vek B Y \,\,\right]
 = \vek{A\wedge B + D_X B-D_Y A}{{[X,Y]}}\,\,,
\ee
where we use $A\wedge B$ to denote  the function of $z$ obtained from $A(z)$ and $B(z)$ by taking the Lie bracket in $\gl$ pointwise, as opposed to $[X,Y]$, which denotes the Lie bracket of $X, Y$ considered as vector fields in $\gl$. The notation $D_A X$ has the usual meaning: the derivative of $X=X(z)$ (at $z$) in the direction of $A=A(z)$.

Let us write $Q=Q(z,z)$ for the vector field in $\gl$ given by the vector field
$q^k(z,z)\frac{\partial}{\partial z^k}$.

For simplicity we will work out the case when $h^{kl}$ is of rank one, which means that the random forcing is applied only in one direction, which will be denoted by $F$ (and considered as a constant vector field in $\gl$). Hence
\be\la{F}
h^{kl}=F^kF^l\,.
\ee
In this case the vector fields for the H\"ormander condition calculation can be taken as
\be\la{21}
\vek 0 F \,\,,\quad \mbox{and}\quad \XX_0=\pul\vek z {Q-\nu z}\,\,,
\ee
We have
\be\la{22}
\left[\vek 0 F, \vek{z}{Q-\nu z}\right]=\vek F {D_F Q - \nu F}\,\,
\ee
and
\be\la{23}
\left[\vek 0 F\,,\left[ \vek 0 F \,, \,\vek z {Q-\nu z}\right]\right]=\vek 0 {D^2_F Q}\,\,.
\ee
This means that we have extended our list of vector fields by the field
\be\la{24}
\vek 0 G \,\,,\,\,G=\pul D^2 Q= Q(F, F).
\ee
We can now take
\be\la{25}
\left[\vek 0 G\,,\left[ \vek 0 F \,, \,\vek z {Q-\nu z}\right]\right]=\vek 0 {D_GD_F Q}\,\,.
\ee
and extend our list of fields by
\be\la{26}
\vek 0 {Q(F,G)}\,.
\ee
We note that the new fields obtained in this way are ``constant" (in the coordinates we use), so the procedure can be easily iterated.

 \begin{definition}
 We will say that $Q$ is  non-degenerate with respect to a set $\mathcal F\subset \gl$ if there is no non-trivial subspace $M\subset \gl$  contaning $\mathcal F$ which is invariant under $Q$, in the sense that $Q(z,z')\in M$ whenever $z,z'\in M$.
\end{definition}

\noindent
We can now formulate the main result of this subsection:

\begin{theorem}\label{thm1}
{The operator of the Fokker-Planck equation~\rf{fpd} satisfies the H\"ormander condition if and only if $Q$ is non-degenerate with respect to the range of the matrix $h$ (considered as a map from $\gl$ to $\gl$).}
\end{theorem}

\noindent
{\it Proof:} The necessity of the condition can be seen when we consider functions depending only on $z$. If there is a non-trivial linear subspace invariant under both $Q$ and the diffusion, then particle trajectories starting at $M$ clearly cannot leave $M$, and therefore the operator cannot satisfy the H\"ormander condition.

On the other hand, if $Q$ is non-degenerate with respect to the range of $h$, then the above calculation shows that the Lie brackets of the fields~\rf{21} (with perhaps several fields of the same form as the first one) generate the fields of the form
\be\la{27}
\vek 0 {X_j}\,,\quad j=1,\dots n\,,
\ee
where $X_1,\dots, X_n\in \gl$ form a basis of $\gl$.
Formula~\rf{21} now shows that the fields of the form
\be\la{28}
\vek {X_j} {Y_j(z)}
\ee
can also be generated. Together with the fields~\rf{27} they clearly form a basis of $T(TG)$ at each point $(a,z)$, and the proof is finished.

\bb
\noindent
{\it Remarks\,\,} 1. If one replaces the damping term $-\nu z$ in~\rf{dl3} by a more general expression $-\nu Dz$, where $D$ is a positive-definite matrix, interesting new questions arise. We plan to address these in a future work.\\
2. Very recently we learned about the paper~\cite{HerzogMattingly}. The results there could be used (modulo simple modifications) to prove the above theorem and also to say more about the set where the solutions of the Fokker-Planck equation are positive.

\begin{corollary}
When $G$ is compact and $Q$ is non-degenerate with respect to the diffucion matrix $h$ in the sense above, the process~\rf{dl3} (and hence also~\rf{DL2}) is ergodic with respect to a  distribution density  given by a function which is independent of $a$. In other words, the Lagrangian positions of the ``particles" are uniformly distributed (with respect to the Haar measure) in the limit of infinite time.
\end{corollary}

In our situation this is not hard to prove once the H\"ormander condition is established by following methods in~\cite{Hairer-ergodic, Hairer-Hormander} and~\cite{Khasminskii}.

\bigskip
\noindent
{\it Remark.}
One should be also able to prove convergence to the equilibrium measure following the methods of~\cite{Villani}, but we will not pursue this direction here. It is perhaps worth reminding that in general there is a difference between uniqueness of the ergodic measure and the convergence to equilibrium. A simple example in our context here is provided by the equation
\be\la{ce1}
f_t+f_{x_1}=\pul f_{x_2 x_2}\,.
\ee
considered in the 2d torus.
Note that this equation does not satisfy the parabolic H\"ormander condition, while its spatial part satisfies the elliptic H\"ormander condition.

\subsection{Constrained diffusion in the momentum space}\label{constr33}
The Euler-Arnold equation~\rf{EA} leaves invariant the co-adjoint orbits~\rf{orb} and also the energy levels $\{H=E\}$.  It is therefore of interest to consider perturbations by noise which ``respects" some of the constraints. For example, one can add noise respecting the coadjoint orbit, but not the energy levels. An example of this situation (in the presence on non-holonomic constraints) is considered in~\cite{Ratiu}. It is closely related to stochastic processes on co-adjoint orbits introduced by Bismut~\cite{Bismut}. One can also consider noise which preserve energy levels but not necessarily the co-adjoint orbits, or one can consider noise which preserves both the coadjoint orbits and the energy levels.

We wish to include all these situations in our considerations, and therefore we will consider the following scenario. We assume that we are given a Hamiltonian $H=H(y)$ and a manifold $M\subset\gls$ which is invariant under the evolution by the Euler-Arnold equation~\rf{EA}. If $M$ is given locally as a non-degenerate level set of some conserved quantities $\phi_1,\dots, \phi_r$ (in the sense that $\{H,\phi_k\}=0\,,\,\,k=1,\dots,r\,$),
\be\la{4-1}
M=\{y\in\gls\,\,,\,\,\phi_1(y)=c_1,\phi_2(y)=c_2,\dots \phi_r(y)=c_r\}\qquad \hbox{locally in $y$},
\ee
there is a natural  measure $m$ on $M$ (invariant for the Hamiltonian flow) which is generated by the volume in $\gls$ (given by our metric there) and the conserved quantities by first restricting the volume measure in $\gls$ to
\be\la{4-2}
M_\ve=\{y\in\gls\,,\,\phi_1(y)\in(c_1-\ve,c_1+\ve),\phi_2(y)\in(c_2-\ve,c_2+\ve),\dots, \phi_r(y)\in(c_r-\ve,c_r+\ve)\}
\ee
then normalizing the restricted measure by a factor $\frac 1{2\ve}$ and taking the limit $\ve\to 0_+$. In the case $r=1$ we have
\be\la{4-3}
m= \frac{1}{|\nabla\phi_1|}\,\,\HH^{n-1}|_M\,,
\ee
where $\HH^{n-1}$ is the $n-1$ dimesnional Hausdorff measure generated by our metric, and the gradient and its norm in the formula are also calculated with our given metric. For general $r$ we have similar formulae, the corresponding expression can be seen easily from the co-area formula, for example. However, the above definition via the limit $\ve\to0_+$ is perhaps more natural, as is relies only on the objects which are ``intrinsic" from the point of view of the definition of $m$: the underlying measure in $\gls$ and the constraits $\phi_k$. (The proof that the limit as $\ve\to 0_+$ is well-defined is standard and is left to the interested reader.)

As the Hamiltonian evolution in the phase-space $T^*G\sim G\times \gls$ preserves the Liouville measure, which, in the $(a,y)$ coordinates defined by~\rf{1}, is the product of the Haar measure on $G$ and the canonical volume measure in $\gls$,
we see that the product of the (left) Haar measure $h_G$ on $G$ and $m$  is an invariant measure for the Hamiltomian evolution in the subset of $T^*G$ given by $G\times M$ in the $(a,y)$ coordinates. If the group $G$ is not unimodular\footnote{Recall that a group is unimodular of the notions of left invariant and right invariant Haar measures coincide. This is the same as demanding that the maps $y\to \Ad a^* y$ preserve the volume in $\gls$, i.\ e.\ have determinant $1$.}, the measure $m$  may not be preserved by the Euler-Arnold equation~\rf{EA} in $\gls$, which represents the symplectic reduction of the original full system. This is because while the
 vector field
 \be\la{4-4}
 v^ke_k+q_k\frac{\partial}{\partial y_k}
 \ee
 in the Liouville equation~\rf{Liouville} is div-free in $T^*G$, its two parts may not be div-free in $G$ or $\gls$ respectvely, unless the group is unimodular.

 The Liouville equation for the evolution is $G\times M$ is the same as~\rf{Liouville}
 \be\la{LiouvilleM}
 f_t+ v^ke_kf+b_k\frac{\partial f}{\partial y_k}=0\,,
 \ee
 where $f=f(a,y)$ now denotes the density with respect to the measure $h_G\times m$ (where $h_G$ is again the left Haar measure on $G$).

 We now consider stochastic perturbations of the Liouville equation~\rf{LiouvilleM} on $G\times M$. As in the Langevin-type equations, the random forces will act only in the $y-$component, so that the kinematic equation $(a^{-1}\dot a)^k=v^k$ is left unchanged.

 We will demand that the stochastic term will also leave invariant the measure $h_G\times m$, and as it acts only in the $y-$variable, it then must leave invariant the measure $m$.

 There is more than one way in which noise can be introduced in a reasonable way into~\rf{LiouvilleM}. For example, if $V$ is a vector field (with coordinates $V_k$) tangent to $M$ which generates a flux on $M$ preserving the measure $m$, one can replace the equation~\rf{H} by
 \be\la{H2}
 d y_k = \{H,y_k\}\,dt+\ve \,V_k\circ dW\,,
 \ee
 where $W$ is the standard $1d$ Wiener process and $\circ$ indicates, as usual, that the corresponding stochastic integrals should be taken in the Stratonovich sense.\footnote{Note that with It\^o integration it the particle trajectories might not stay in $M$.}
 The corresponding Fokker-Planck equation is given by
 \be\la{fpm1}
 f_t+ v^ke_kf+b_k\frac{\partial f}{\partial y_k}=\frac{\ve^2} 2(D_V^*)^2f\,,
 \ee
 where $D_V^*$ is adjoint to $D_V=V_k\frac{\partial}{\partial y_k} $  with respect to the measure $m$.
 As the flux by $V$ preserves $m$, we in fact have $D_V^*=-D_V$. In this case the operators $D_V^2$ and $(D_V^*)^2$ coincide and arise from the functional
 \be\la{funcV}
 \int_M \pul |D_V f|^2\, d\,m
 \ee
 This is in some sense the ``minimal non-trivial noise" model, and it might be of interest in some situations.

 Here we will consider the situation when the noise is non-degenerate in $M$, leaving the interesting case of the degenerate noise in $M$  to future studies. Our motivation is the following. For $\ve>0$ we consider the usual Brownian motion in $\gls$, but restricted to the set $M_\ve$ above, with the understanding that the trajectories ``reflect back" (we can think about an action of some control mechanism) at the boundary (corresponding to the Neumann condition at the boundary for the corresponding Fokker-Planck equation, which is just the heat equation in this case).   As $\ve\to 0_+$, a good model for the limiting process on $M$ is given by an operator constructed as follows.
 First, we take the metric induced on $M$ by the given metric in $\gls$. Assume the metric is given by $\tilde g_{ij}$ in some local coordinates. Assume the measure $m$ is given as $m(x)\,dx$ in these coordinates, where $m(x)$ is a (smooth) function. Denoting  by $\tilde g$ the determinant of $\tilde g_{ij}$, the volume element given by $\tilde g_{ij}$ in our coordinates is $\sqrt{\tilde g}\,dx$. We then define a new metric
 \be\la{newh}
 h_{ij}=\varkappa \tilde g_{ij}
 \ee
 so that the volume element $\sqrt{h}\,dx$ satisfies
 \be\la{hdef}
 \sqrt{h}\,dx=m(x)\,dx\,\,.
 \ee
 Then we take the generator of our process to be the Laplace operator on $M$ with respect to the metric $h$. We will denote this operator by $L_M$.
 Our Fokker-Planck equation then will be
 \be\la{fpm2}
 f_t+ v^ke_kf+b_k\frac{\partial f}{\partial y_k}=\frac{\ve^2}2 L_M f\,.
 \ee
 We will be interested in ergodicity properties of the process given by this equation.

\bb
In the remainder of this subsection we will assume again~\rf{HH}, i.\ e.\ the Hamiltonian $H$ is quadratic (and positive definite). We can then ``lower the indices" and work with $TG$ and $\gl$ rather then with $T^*G$ and $\gls$. We will denote by $Z$ the image of $M$ in $\gl$ under the ``lowering indices" map, and will denote the elements of $Z\subset\gl$ by $z$, with coordinates $z^k$. Similarly to~\rf{vey} we have
$z^k=v^k$.
The Fokker-Planck equation~\rf{fpm2}, now considered on $G\times Z$ becomes

\be\la{fpz}
f_t+z^ke_k f + q^k(z,z)\frac{\partial f}{\partial z^k} = \frac{\ve^2}{2}Lf \,\,,
\ee
where $q^k$ is defined by~\rf{symq}, and $L$ is the operator on $Z$ corresponding to $L_M$. It is of course again a Laplacian for some metric on $Z$ (which is conformally equivalent to the metric on $Z$ induced by the underlying metric in $\gl$).

Let us now consider conditions under which the operators corresponding to \rf{fpz} or \rf{fpm2} satisfy the usual H\"ormander commutator condition for hypoellipticity.

\begin{definition}\label{def1}
 A p-hull\,\,\footnote{Here p stands for parabolic, as the definition is tied to the parabolic H\"ormander condition.} of a subset $S\subset \gl$ is the smallest Lie sub-algebra $\mathfrak h\subset\gl$  with the following properties:
 \begin{itemize}
  \item[(i)] $\mathfrak h$ contains the set $S-S=\{s_1-s_2, \,,\,s_1,s_2\in S\}$\,,
 \item [(ii)] $\mathfrak h$ is invariant under the mappings $\Ad s\colon  \,\,z\to [s,z]$ for each $s\in S$.
 \end{itemize}
\end{definition}

\noindent
Remarks:

\sm
\noindent
1. The p-hull will be relevant in the context of the evolution equation~\rf{fpz}. For the ``spatial part" of the operator~\rf{fpz}, obtained by omitting the term $f_t$, the relevant ``hull" is simply the Lie algebra generated by $S$.

\sm
\noindent
2. Condition (i) in the definition already implies that $\mathfrak h$ is invariant under the mapping $\Ad (s_1-s_2)$ for any $s_1,s_2\in S$. Therefore in (ii) we can require invariance of $\mathfrak h$ under $\Ad s_0$ for just one fixed $s_0\in S$ (and - given (i) - the definition will be independent of the choice of $s_0$).

The main result of this section is the following:

\begin{theorem}\label{thm2}
In the notation introduced above, assume that $M$ is a smooth analytic submanifold of $\gls$. Then the operator on $G\times M$ corresponding to~\rf{fpm2} (or, equivalently, the operator on $G\times Z$ corresponding to~\rf{fpz}) satisfies the H\"ormander condition if and only if the p-hull of $Z$ coincides with $\gl$.
\end{theorem}

\noindent
{\it Proof:}

\noindent
Let us first show that the p-hull condition is necessary for the H\"ormander condition. One can see this from the Lie bracket calculations below, but it is instructive to verify it directly.
Assume $\mathfrak h\,$ is a non-trivial Lie subalgebra of $\gl$ containing $Z-Z$ for which we can find $z_0\in Z\setminus\hl$  such that $\hl$ is invariant under $\Ad z_0$. Let us set
\be\la{edef}
e=z_0^ke_k\,\,.
\ee
The Lie algebra $\hl$ defines (locally) a foliation $\mathcal F$ of $G$ into cosets $aH$, where $H$ is the (local) Lie subgroup of $G$ corresponding to $\hl$. The main point now is that the invariance of $\hl$ under $\Ad z_0$ implies that the flow given by the equation
\be\la{eq1}
a^{-1}\dot a =e
\ee
preserves the foliation. (Another formulation of this statement could be that the equation~\rf{eq1} ``descends" to $G/H$.)
This means that the perturbations given by the stochastic terms in~\rf{fpz} will still preserve the foliation (e.\ g.\ by the Stroock-Varadhan Theorem~\cite{StroockVaradhan}) and it is not hard to conclude that set of points reachable by the corresponding process cannot be open.

For the proof that the p-hull condition is sufficient, we write our operator (locally) in the form
\be\la{op}
f_t+ \XX_0 f- \sum_{j=1}^m \XX_j^2 f\,,
\ee
where $m$ is the dimension of $Z$ (which is of course the same as the dimension of $M$) and $\XX_j$ are suitable vector fields on $G\times Z$. All these fields will be of the form~\rf{fields}, and we will use the same notation as in~\rf{fields} in what follows.
We will be working locally near a point $(a,z_0)\in G\times Z$.
We choose $\XX_j\,,j=1,\dots, m$ so that
\be\la{f1}
\XX_j=\vek 0 {Y_j}\,,\qquad j=1,\dots, m\,.
\ee
where $Y_j$ are analytic near $z_0$ and  $Y_j(z)$ form a basis of $T_z Z$ for each $z$ close to $z_0$. The field $\XX_0$ will be of the form
\be\la{f2}
\XX_0=\vek z V\,,
\ee
where $V$ is an field on $Z$ (analytic near $z_0$).
Let us consider local analytic vector fields on $G\times Z$ near $(a,z_0)$ of the form
\be\la{ourfields}
\XX(a,z)=\vek {X(z)}{Y(z)}
\ee
as a (Lie) module $\AAA$ over the set of analytic functions of $z\in Z$. (Recall that we assume that $Z$ is analytic.)
Let $\MM$ be the minimal (Lie) submodule of $\AAA$ satisfying the following requirements:
\begin{itemize}
\item[(a)] $\MM$ contains $\XX_1,\XX_2,\dots, \XX_m$\,, and
\item[(b)]  $\MM$ is invariant under the map $\Ad \XX_0\colon \XX\to [\XX_0, \XX]$\,, where $[\,\cdot\,,\,\cdot\,]$ denotes the Lie bracket for vector fields.
\end{itemize}
The parabolic H\"ormander condition at $(a,z_0)$ for the fields $\XX_0,\XX_1,\dots, \XX_m$ then is that
 \be\la{pH}
 \{\XX(a,z_0)\,,\,\,\XX\in \MM\} = T_{(a,z_0)}( G\times Z)\,\,.
\ee
For $\XX\in \AAA$ we will denote by $\pi \XX\in\gl$ the projection to the first
component in the notation~\rf{fields}, i.\ e.\
\be\la{pidef}
\pi\vek X Y = X\,.
\ee
Let
\be\la{Mdef}
M=\pi \MM\,.
\ee
As $\MM$ contains the vector fields~\rf{f1}, the condition~\rf{pH} is equivalent to
\be\la{pH1}
M_{z_0}=\{X(z_0)\,,\, X\in M\}= \gl\,.
\ee
Using~\rf{lbc} and the fact that the fields $\XX_1,\dots, \XX_m$ belong to $\MM$,  we see that $M$ has the following properties.

\be\la{p0}
\hbox{\sl If $Y$ is an analytic vector field on $Z$ (defined locally near $z_0$), then $Y\in M$.}
\ee
This follows by taking the Lie bracket of $\vek 0 Y$ and $\XX_0$.
\be\la{p1}
\begin{array}{c}
\hbox{\sl If $A\in M$ and $Y$ is an analytic vector field on $Z$ (defined locally near $z_0$),} \\ \hbox{\sl then $D_Y A$ is in $M$.}
\end{array}
\ee
This follows by taking the Lie bracket of $\vek 0 Y$ and  an $\XX$ with $\pi \XX=A$.

\be\la{p2}
 \hbox{\sl If $A,B\in M$, then $A\wedge B$ is in $M$.}
 \ee
  This follows by taking the Lie bracket of the fields $\XX$ and $\YY$ with $\pi \XX=A$ and $\pi\YY=B$ and then using~\rf{p1}.

\be\la{p3}
\hbox{\sl If $A\in M$, then $z\wedge A$ is in $M$.}
 \ee
 This follows by taking the Lie bracket of $\XX$ with $\pi \XX= A$ with $\XX_0$ and using~\rf{p0}, \rf{p1}, and \rf{p2}.
Taking these properties of $M$ into account, it is clear that the proof of the theorem will be finished if we show that
\be\la{ZmZ}
Z-Z\subset M_{z_0}\,.
\ee
Let $l$ be a linear function in $\gl$ which vanishes on $M_{z_0}$. As $Z$ is analytic, the function $l$ considered as a function on the manifold $Z$ will be analytic. The property~\rf{p1} of $M$ implies that the derivatives of all orders $\ge 1$  of $l$ at $z_0$ vanish, and therefore $l$ must be constant on $Z$. In particular, $l$ must vanish on $Z-Z$. We see that no point of $Z-Z$ can be separated from the subspace $M_{z_0}$ by a linear function, and~\rf{ZmZ} follows. This finishes the proof of the theorem.

\begin{corollary}
If the assumptions of Theorem~\ref{thm2} are satisfied and the group $G$ is compact, then any solution of the Fokker-Planck equation~\rf{fpm2} approaches constant. In particular, the system is ergodic  for the (stochastic) dynamics, with the unique ergodic measure given by the constant density $f$.
\end{corollary}

\noindent
{\it Proof:}
We note that
\be\la{decay1}
\frac{d}{d\, t}\int_{G\times Z} f^2(a,z,t) \,da\,m(dz)=-\ve^2\int_{G\times Z}|\nabla_z f(a,z,t)|^2\,da\,m(dz)\,,
\ee
where we take $Z$ with the metric defining $L$. By regularity which follows from the H\"ormander condition we can consider the $\Omega-$limit set $\Omega(f_0)$ of the evolution starting with $f_0$, and it consists of smooth functions. Moreover the integral on the right of~\rf{decay1} has to vanish identically for each function in $\Om(f_0)$, by the usual Laypunov-function-type arguments. This means that any function in $\Omega(f_0)$ is constant in $z$ and hence solves the equation
\be\la{decay2}
f_t+z^ke_kf=0\,.
\ee
It is now easy to see that our assumptions imply that such $f$ is constant also in $a$.

\subsection{A calculation for a non-compact group}
We now consider the situation in the previous subsection for the special case $G=\R^n$ and a one-dimensional manifold $Z\subset \gl\sim\R^n$. In other words, $Z$ will be an analytic curve in $\R^n$. We will see that analyticity is not really needed for the calculation below, but we keep it as an assumption, so that we have the H\"ormander condition for the Fokker-Planck equation under the assumptions of Theorem~\ref{thm2}. We will assume that $Z$ is equipped with a measure $m$, the density of which is also analytic with respect to the parameter which gives an analytic parametrization of $Z$. We will re-parametrize $Z$ so that it is given by an analytic periodic function
\be\la{ nc1}
\gamma\colon R \to Z\subset \R^n
\ee
with minimal period $l$ and, in addition, the measure (as measured by $m$)  of a segment on the curve between $\gamma(s_1)$ and $\gamma(s_2)$ for some $0\le s_1<s_2<l$ will be given  by $s_2-s_1$. Sometimes we will also write
\be\la{nc2}
\gamma(s)=z(s)\,,
\ee
with slight abuse of notation which will hopefully not cause any confusion.
In this special case the Fokker-Plack equation discussed in the previous section, written in the variables $a=(a^1,\dots, a^n)\in G$ and $s$ (which parametrizes $Z$), is \be\la{fps}
f_t+z^k(s)\frac{\partial f}{\partial a^k} = \frac{\ve^2}2 \frac{\partial ^2 f}{\partial s^2}\,,
\ee
where $f=f(a^1,\dots,a^n,s,t)\,$ is periodic in $s$, with period $l$.
The p-hull condition from Definition~\ref{def1} is that $Z-Z$ generates $\R^n$.

We are interested in the long-time behavior of the solutions of~\rf{fps}. We will assume that the p-hull condition is satisfied. It is easy to see that the case when the condition is not satisfied can be reduced to this case by a suitable choice of variable.\footnote{Here and below this is of course meant only in the context of the example we are considering in this subsection. }

We note that the change of variables $a^k\to a^k - z^k_0t$ for some $z_0\in\R^n$ is equivalent to shifting $Z$ to $Z-z_0$. We can therefore assume without loss of generality that
\be\la{nc3}
\int_0^l \gamma(s)\,ds = \int_Z z\,m(dz)=0\,.
\ee
This condition enables us to write
\be\la{phidef}
\gamma(s)=\vf''(s)
\ee
for some periodic (analytic) $\vf\colon \R\to\R^n$.
An important role will be played by the matrix
\be\la{covm}
\Sigma_{kl}=\frac 1 l\int_0^l \vf'_k(s)\vf'_l(s)\,ds\,.
\ee
\vbox{
\begin{proposition}\label{prop1}
Assume~\rf{nc3} (which can be always achieved by a change of variables $a\to a-z_0t$) and let $\Sigma_{kl}$ be defined by~\rf{covm}. For any compactly supported initial density $f_0=f_0(a,s)$ (normalized to total mass one) the quantity
\be\la{limf}
a\to {t^{\frac n2}} \int_0^lf(\sqrt{t}\, a, s,t)\,ds
\ee
converges as $t\to\infty$ (in distribution) to the density of the normal distribution with average $0$ and covariance matrix $\frac{4}{\ve^2}\Sigma_{kl}$.
In other words, the distribution of the positions of trajectories starting at time t in some compact region will approach (after re-scaling) the same distribution as the diffusion with covariance matrix $\frac{4}{\ve^2}\Sigma_{kl}$.
\end{proposition}
}
\noindent{\it Proof:}

\noindent
We will work with the corresponding stochastic ODE
\be\la{sODE}
\begin{array}{rcl}
\dot a & = & \gamma(s)\\
\dot s & = &  \ve \dot w\,,
\end{array}
\ee
where $w(t)$ is the standard one-dimensional Wiener process starting at the origin.
Our task reduces to evaluating
\be\la{integral}
a(t)-a(0)=\int_0^t \gamma(\ve w(t'))\,dt'=\int_0^t\vf''(\ve w(t'))\,dt'\,.
\ee
We will evaluate the integral by a standard procedure based on the martingale version of the central limit theorem. We only sketch the main steps.
By It\^o formula we have
\be\la{nc4}
\vf(\ve w(t))-\vf(\ve w(0))=\int_0^t \ve\vf'(\ve w(t'))dw(t') +\int_0^t\frac{\ve^2}2\vf''(\ve w(t'))\,dt'\,.
\ee
We  re-write this as
\be\la{nc5}
\frac1{\sqrt{t}} \int_0^t \gamma(\ve w(t'))\,dt'=\int_0^t\frac {2}{\ve\sqrt {t}}\vf'(w(t'))(-dw(t'))+\frac2{\ve^2\sqrt{t}} \left(\vf(\ve w(t))-\vf(\ve w(0))\right)
\ee
The last term on the right clearly approaches zero for $t\to\infty$, as $\vf$ is bounded.
The key point now is to use a martingale version of the central limit theorem (such as, for example Theorem 3.2, page 58 in~\cite{Hall}) to get a good asymptotics for the  integral on the right. The covariance matrix for that integral generated along a trajectory $w(t')$ is
\be\la{var}
\frac {4}{\ve^2 t} \int_0^t \vf'_k(\ve w(t'))\vf'_l(\ve w(t'))\,dt'\,.
\ee
For large times $t'$ the distribution of the variable $\ve w(t')$ taken mod $l$ will be approaching the uniform distribution in $[0,l)$  and therefore it is not hard to see that for the purposes of our calculation we can replace the random quantity~\rf{var} by a deterministic quantity given by
\be\la{var2}
\frac {4}{\ve^2l}\int_0^l \vf'_k(s)\vf'_l(s)\,ds= \frac{4}{\ve^2}\Sigma_{kl}\,.
\ee
The claim of the proposition now essentially follows from the central limit theorem.

\bigskip
\bigskip
\centerline{\bf Acknowledgement}

\smallskip
\noindent
We thank Jonathan Mattingly for an illuminating discussion.\\
The research was supported in part by grants DMS 1362467 and DMS 1159376 from the National Science Foundation.


\bigskip
\begin{thebibliography}
{99}

\bibitem{Arnold}  Arnold, V.\ I., {\sl  Sur la g\'eom\'etrie diff\'erentielle des groupes de Lie de dimension infinie et ses applications \`a l'hydrodynamique des fluides parfaits},  Ann.\ Inst.\ Fourier (Grenoble) 16 1966 fasc. 1, 319--361.
\bibitem{ArnoldKhesin}
Arnold, V.\ I., Khesin, B.\ A., {\sl Topological Methods in Hydrodynamics,} Springer, 1998.
\bibitem{Bismut} Bismut, J.-M., {\sl Mecanique Al\'eatoire,} Lecture Notes in Math., vol.\ 866, Springer 1981.
\bibitem{Einstein1905}
Einstein, A., {\sl \"Uber die von der molekularkinetischen Theorie der W\"arme geforderte Bewegung von in ruhenden Fl\"ussigkeiten suspendierten Teilchen,} Annalen der Physik, 322 (8), 549--560.
\bibitem{Jurdjevic}
Jurdjevic, V., {\sl  Geometric Control Theory}, Cambridge Studies in Advanced Mathematics vol.\ 52, 1997
\bibitem{Hairer-Hormander}
     Hairer, M., {\sl On Malliavin's proof of H\"ormander's theorem,} Bull.\ Sci.\ Math.\ 135 (2011), no. 6-7, 650--666
\bibitem{Hairer-ergodic} Hairer, M., {\sl Ergodic theory for stochastic PDEs,} Lecture notes from the LAM-EPSRC short course held in July 2008.

\bibitem{HM}
Hairer, M., Mattingly, J.\ C.,\,  {\sl Ergodicity of the 2D Navier-Stokes equations with degenerate stochastic forcing}, Ann. of Math. (2) 164 (2006), no. 3, 993--1032.

\bibitem{Hall}
Hall, P., Heyde, C.C.,{\sl  Martingale Limit Theory and its Application,} Academic Press,
1980.

\bibitem{HerzogMattingly}
Herzog, D.\ P., Mattingly, J.\ C.,\,
\,{\sl A practical criterion for positivity of transition densities,} Nonlinearity, 28 (2015) 2823--2845.
\bibitem{Ratiu}
Hochgerner, S., Ratiu, T., {\sl Geometry of non-holonomic diffusion,} arXiv:1204.6438.

\bibitem{Hormander}
H\"ormander, L., {\sl Hypoelliptic second order differential equations,} Acta Math.\ 119, 1967, 147--171.
 \bibitem{Khasminskii}
 Khasminskii, R., {\sl Stochastic stability of differential equations,} second edition. Stochastic Modelling and Applied Probability, 66. Springer,  2012.
\bibitem{Kuksin}
 Kuksin, S.\ B.,  {\sl Randomly forced nonlinear PDEs and statistical hydrodynamics in 2 space dimensions. Zurich Lectures in Advanced Mathematics,} European Mathematical Society (EMS), Z\"urich, 2006.
\bibitem{MarsdenWeinstein}
Marsden, J., Weinstein, A., {\sl  Coadjoint orbits, vortices, and Clebsch variables for incompressible fluids,} Order in chaos (Los Alamos, N.M., 1982), Phys. D 7 (1983), no. 1-3, 305--323.
\bibitem{StroockVaradhan}
 Stroock, D.\ W., Varadhan, S. R. S., {\sl On the support of diffusion processes with applications to the strong maximum principle,} Proceedings of the Sixth Berkeley Symposium on Mathematical Statistics and Probability (Univ. California, Berkeley, Calif., 1970/1971), Vol. III: Probability theory, pp. 333--359.
\bibitem{Villani}
Villani, C., {\sl Hypocoercivity,}
Mem.\ Amer.\ Math.\ Soc.\ 202 (2009), no.\ 950.


\end{thebibliography}
\end{document}